\magnification1100

\input amstex
\documentstyle{amsppt}
\topmatter
\title{On increasing subsequences of i.i.d. samples}
\endtitle
\author{Jean-Dominique Deuschel and Ofer Zeitouni}
\endauthor
\affil{Technische Universit\"at
Berlin and Technion}\endaffil
\address{ Fachbereich Mathematik, TU-Berlin,
 Stra\b e des 17. Juni 135, D-10623 Berlin}\endaddress
\address{ Department of Electrical Engineering, Technion,
Haifa 32000, Israel}\endaddress
\thanks{ The second author acknowledges support from
a US-Israel B.S.F. grant}\endthanks
\email{deuschel@@stoch1.math.tu-berlin.de}
\endemail
\thanks{{\it AMS 1991 subject classifications}. Primary 60G70; secondary
60F10.}\endthanks
\thanks{{\it Key words and phrases.} 
Increasing subsequences, random permutations,
large deviations.}\endthanks
 \email{zeitouni@@ee.technion.ac.il}\endemail
\abstract{We study the fluctuations, in the large deviations regime,
of the longest 
increasing subsequence of a random i.i.d. sample on the unit square.
 In particular,
our results yield the precise upper and lower exponential tails for the length of
the longest increasing subsequence of a random permutation. }\endabstract
\endtopmatter

\heading{\S. 1 Introduction}\endheading
Let $\{Z_i\}_{i=1}^n =\{(X_i,Y_i)\}_{i=1}^n$ denote a sequence of i.i.d.
random variables
with marginal law $\mu$ on the unit square $\bold Q=[0,1]^2$.
Throughout, we make the assumption that $\mu$ possesses
a strictly positive density  $p\in C^1(\bold Q)$ with respect to the
Lebesgue measure $\lambda$ on $\bold Q$.

A subsequence $\{Z_{i_1},...,Z_{i_\ell}\}\subseteq \{Z_i\}_{i=1}^n$ is
called a {\it monotone increasing
subsequence of length} $\ell$, if
$$X_{i_j}<X_{i_{j+1}}\quad\text{and}
\quad Y_{i_j}<Y_{i_{j+1}}, \quad\text{for}\quad j=1,...,\ell-1.$$
Define next $\ell_{\max}(n)$ to be the length of the longest
increasing subsequence in the
sample $\{Z_i\}_{i=1}^n$. {Note that we do not
require that $i_j<i_{j+1}$.

In the case that $\mu=\lambda$,  $\ell_{\max}(n)$ possesses the
same law as the length of the longest
increasing subsequence of a random permutation,
denoted hereafter by $L_{\max}(n)$. Building on the fact
that
$$\lim_{n\to\infty} \frac{L_{\max}(n)}{\sqrt n} =
2\qquad\text{in probability},$$
c.f. \cite{10},\cite{11},
we showed in \cite{3} that
$$\lim_{n\to\infty}
\frac{\ell_{\max}(n)}{\sqrt n}
 = 2\bar J_\mu\qquad\text{in probability},\tag{1.1}$$
where $\bar J_\mu\in\Bbb R^+$ is the solution to
 the variational problem
$$\bar J_\mu=\sup_{\phi\in B^\uparrow}
\int_0^1\sqrt{p(x,\phi(x))\dot\phi(x)}\,dx\tag{1.2}$$
with
$$B^\uparrow\equiv\{\phi:[0,1]\longrightarrow[0,1],
\quad\text{non-decreasing, absolutely continuous}\}.$$
Furthermore, it follows from \cite{3} that any
longest increasing subsequence will
concentrate on the solutions to the variational problem (1.2).
See \cite{4}, Proposition 4.4, for an alternative expression for $\bar J_\mu$.

Note that $\bar J_\mu=1$ for $p(x,y)=1$,
in which case 
the maximum is achieved on the diagonal $x\longrightarrow
\phi(x)=x.$

The fluctuations of $\ell_{\max}(n)$ and $L_{\max}(n)$ are
highly nontrivial and have been investigated in
several papers, c.f. \cite{2}, \cite{9}, \cite{1}, \cite{5}. In particular,
Aldous and Diaconis have exhibited quite different
behaviors in their upper and lower tails.
Our goal in this paper is to provide information on {\sl the large
deviations of these fluctuations}. The results and
techniques differ sharply in the study of lower and
upper tail, and we divide the discussion in the rest of this
introduction between these two cases.

Turning our attention to {\sl the lower tail}
 we first show in Theorem 1 that for any $-2<c<0$,
$$\lim_{n\to\infty}\frac1n\log P( L_{\max}(n)<(2+c)\sqrt{n})=-2H_0(c ),
\tag{1.3}$$
with an explicit function $H_0$, 
first introduced by
Logan and Shepp in \cite{6},
$$H_0(c )=-\frac12 +\frac{(2+c)^2}{8}+
\log\frac{c+2}{2}-\big(1+\frac{(c+2)^2}{4}\big)
\log\big(\frac{2(c+2)^2}{4+(c+2)^2}\big)\,.$$
See Fig. 1 for a plot of $H_0(\cdot)$.
The proof based on the random Young tableau correspondence is
purely combinatoric and sheds no light on 
the random mechanism responsible for the large deviations.
In particular, 
we cannot prove that, conditioned on $L_{\max}(n)<(2+c)\sqrt{n}$, the
longest increasing subsequence concentrates around a curve $\phi(\cdot)$.

While we could hope to 
use this result in order to prove an exponential
lower tail for general $\mu$,
that is
$$\lim_{n\to\infty}\frac1n
\log P( \ell_{\max}(n)<(2\bar J_\mu+c)\sqrt{n})=-2H_\mu(c ),\tag{1.4}$$
with $H_\mu( c)>0$ for $-2\bar J_\mu<c<0$, we were not able to compute $H_\mu$
explicitly, nor to prove the existence of
the limit in (1.4). We thus present in Propositions 2.1 and 2.2  nontrivial
upper and lower bounds on the left hand side of (1.4), avoiding
the question of existence of the limit.

\bigpagebreak

The situation is quite different for {\sl the upper tail}, here an
easy sub-additive argument shows that
for $c>0$
$$\lim_{n\to\infty}\frac1{\sqrt n}\log P(L_{\max}(n)>(2+c)\sqrt n)=-U_0(c ),
\tag{1.5}$$
for some nontrivial convex rate function $U_0(\cdot)$.
In the first version of this work, we presented only
bounds on $U_0(\cdot)$, leaving open the explicit evaluation of
this function.
Subsequently, T. Sepp\"{a}l\"{a}inen has proved in \cite{8}, using
Hammersley's particle system associated with the 
Poissonized version of $L_{\max}(n)$, that
$$ U_0(c)=\beta(c):= 2(2+c)\cosh^{-1}(c/2+1)-4
\sqrt{c^2+4c}\,.$$
See Fig. 2 for a plot of $U_0(\cdot)$.
In fact, Kim \cite{5} had already observed, by combinatorial techniques,
the upper bound in (1.5) with the
function $U_0(c)$. For the sake of completeness, we will present in
Section 3 a 
combinatorial proof of the lower bound in (1.5).

Our interest is in exploring the similar
 question for $\ell_{\max}(n)$, where the sub-additive
argument is not applicable. Our main result in this direction 
(c.f. Theorem 3) is in fact that
$$\lim_{n\to\infty}\frac1{\sqrt n}
\log P(\ell_{\max}(n)>(2\bar J_\mu+c)\sqrt n)=-U_\mu(c)\tag{1.6}$$
where
$$U_\mu(c)=\bar J_\mu U_0(\frac{c}{\bar
J_\mu}) .$$
Moreover we show that, under the conditioning that
$$\{\ell_{\max}(n)>(2\bar J_\mu + c)\sqrt n\},$$
the longest increasing subsequences concentrate near the
maximizing curves in (1.2).
 
The precise statement and proof of these theorems is
followed in section 4  by a discussion
and several conjectures and open questions.
\bigpagebreak

{\bf Acknowledgment} We thank JH
Kim for pointing out to us that $H_0(c )$ is an upper
bound in (1.3), A. Dembo for useful discussions, and
T. Sepp\"{a}l\"{a}inen for sending us a copy
of \cite{8}, where a probabilistic proof of (1.5) was
first provided.

\heading{\S. 2 The lower tail}\endheading
In this section we first describe the large deviations for the
uniform measure $p(x,y)=1$.
Let us recall some notation  from \cite{6}: For  $f\in\Cal F$ ,
the class of nonnegative,
nondecreasing functions on $[0,\infty)$ of unit integral, define
$$H(f)\equiv \int_0^\infty
\int_0^{f(x)}\log(f(x)-y+f^{-1}(y)-x)\,dy\,dx\,.\tag{2.1}$$
Let, for $c\in(-2,0)$,
$$H_0(c )\equiv\inf\{H(f),\quad f\in\Cal F, f(0)=2+c\}\ + \ \frac12,$$
c.f. \cite{6}.
Then $H(f)\ge 0$, and
$$H_0(c )=-\frac12 +\frac{(2+c)^2}{8}+
\log\frac{c+2}{2}-\big(1+\frac{(c+2)^2}{4}\big)
\log\big(\frac{2(c+2)^2}{4+(c+2)^2}\big).\tag{2.2}$$
Note that $H_0$ is a strictly convex,
monotone decreasing function with minimum $0$
at $c=0$, c.f. Fig. 1 below.
Our first result, which is an
 immediate consequence of \cite{6}, is based on Schensted's
identity:
\proclaim{Theorem 1} For any $-2<c\le 0$,
$$\lim_{n\to\infty}
\frac1n\log P( L_{\max}(n)<(2+c)\sqrt{n})=-2H_0(c ).$$
\endproclaim
\demo{Proof}
\comment
 The upper bound in (2.3)
 being nothing but \cite{6} (1.9), (1.10) and (3.21), we concentrate
the sequel in proving the lower bound
$$\liminf_{n\to\infty}
\frac1n\log P( L_{\max}(n)<c\sqrt{n})\ge-2H_0(c ).\tag{2.3}$$
\endcomment
The basic idea is to use a combinatorial identity of Schensted,
expressing the probability distribution of
$L_{\max}(n)$ in terms of Young tableaux,
c.f. \cite{7} and \cite{6}, \S 1:

A Young shape $\tau$ of size $|\tau|
=n$
is an array of $n$ unit squares $s$, left and bottom justified,
 whose columns' lengths are nonincreasing
from left to right.
The {\it hook length} $\sigma(s)$ of a square $s$ in the shape $\tau$ is just
 the number
of squares in $\tau$
 directly above and to the right of it,
counting each square  exactly once, c.f. Fig 3.
Let $\pi(\tau )=\prod_{s\in\tau}\sigma(s)$ denotes the hook
product, i.e. the product of all hook
lengths in the tableau $\tau$. Then, the
Schensted identity states that
$$P(L_{\max}(n)=k)=
\sum_{\tau: \tau(0)=k}\frac{n !}{(\pi(\tau))^2},\qquad k=1,...,n\,,
\tag{2.3}$$
where the sum is taken over all shapes $\tau$ containing
$n$ squares,
possessing a first column of length
$k$.
In order to estimate $P(L_{\max}(n)\le (2+c)\sqrt n)$ for
fixed $c\in (-2,0)$ it suffices to find
 an optimal shape $\tau_n$ with $\tau_n(0)\le (2+c)\sqrt n$
which maximizes the hook product
$\pi(\tau_n)$.
This is in essence the argument of \cite{6} which yields the upper bound,
c.f.  (1.9), (1.10) and (3.2) there.
We hence concentrate
in the sequel in proving the lower bound
$$\liminf_{n\to\infty}\frac1n\log
 P( L_{\max}(n)<(2+c)\sqrt{n})\ge-2H_0(c ).$$
Our goal is to find for fixed $c\in (-2,0)$  a sequence of shapes
 $\{\tau_n\}$ of  maximal hook product such that
$\lim_{n\to\infty}|\tau_n|/n=1$ and
 $\lim_{n\to\infty}\tau_n(0)/n^{1/2}\le (2+c)$.
Let  $f_0^c\in\Cal F$ be such that
$$H(f_0^c)=\inf\{H(f):\quad f\in\Cal F\qquad\text{with}\qquad f(0)=2+c\}.$$
The curve $f_0^c$ is constructed in \cite{6}, it has the support
$$b_0(c )=\frac1{(2+c)} -\frac{(2+c)}{4}+\sqrt{2+\frac{(2+c)^2}{2}}.$$
Hence, the length of the curve $\{(x,f^c_0(x)),  0\le x\le b_0(c )\}$ is
bounded by some constant $k_c$.

We   construct a particular Young tableau
 out of $f^c_0$. For $i=1,...,[b_0(c)\sqrt n]\equiv i_{\max}$
set $j(i)=[f^c_0(\frac{i}{\sqrt n})\sqrt n]$.
  Note that $j(i)$ is a decreasing sequence, and, because
the length of  $\{(x,f^c_0(x)),  0\le x\le b_0(c )\}$ is bounded,
$$m_n=\sum_{i=1}^{i_{\max}}\sum_{j=1}^{j(i)} 1\ge n - k_c\sqrt n.$$
The sequence $\{(i,j(i)), i=1,...,i_{\max}\}$ defines a
Young tableau $\tau_n$ of size $m_n$.
Moreover for any $y<f_0^c(x)$ such that
$\bar \imath=[x\sqrt n]\le i_{\max}$ and $\bar \jmath=[y\sqrt n]\le j(i)$,
denoting by $\pi_{\bar \imath \bar \jmath}$ 
the hook length of the square  with
indices $(\bar \imath,\bar \jmath)$,
$$\log\pi_{\bar 
\imath \bar \jmath}(\tau_n)\le\log(f_0^c(x)-y + (f_0^c)^{-1}(y)-x)+\log n.$$
Hence, for some constant $C>0$ independent of $n$,
whose value may change from line to line,
$$\align
nH(f_0^c)&=
n\int_0^{b_0( c)}\int_0^{f_0^c(x)}\log\big(f_0^c(x)-y+(f_0^c)^{-1}(y)
-x\big)\, dx \, dy\\
\ge& -\frac{n}{2}\log n +
 \log\pi(\tau_n) + n\sum_{i=1}^{i_{\max}}\big(f_0^c(\frac{i}{\sqrt n})-
f_0^c(\frac{i+1}{\sqrt n})+\frac1{\sqrt n}\big)\int_0^{1/\sqrt n}
\log x \,dx\\
\geq
&-\frac{n}{2}\log n  + \log\pi(\tau_n)-C\sqrt n\log n\,.
\tag{2.4} \endalign
$$
It follows that for any $-2<c<0$,
$$\align  P(L_{\max}(m_n)<(2+c)\sqrt n)&\ge\frac{m_n !}{(\pi(\tau_n))^2}
\ge \frac{m_n!}{n!} n! e^{-2n\log\sqrt n}e^{-C\sqrt n\log n}
 e^{-2n H(f_0^c)}\\
&\ge e^{-C\sqrt n\log n}e^{-2n[H(f_0^c)+\frac12]}.\endalign$$
Finally, for any $\bar c<2$, by rescaling,
$$ P(L_{\max}(n)<\bar c\sqrt n)
\le P(L_{\max}(m_n)<\bar c\sqrt{m_n})=
P(L_{\max}(m_n)<\bar c\sqrt n\frac{\sqrt{m_n}}
{\sqrt n})\,,$$
and the conclusion follows from the continuity of $H(f_0^c)$ in c.\qed
\enddemo

An immediate corollary, which will be useful below, is the following:
\proclaim{Corollary 1} For any $-2<c< 0$ there exists a function
$\eta(c,\delta)$ satisfying 
$$\lim_{\delta\to 0} \eta(c,\delta)=0$$
such that if $p(x,y)$ satisfies $(1-\delta)\leq p(x,y)\leq (1+\delta)$ then
$$\limsup_{n\to\infty}
|\frac1n\log P( \ell_{\max}(n)<(2+c)\sqrt{n})+2H_0(c )|
\leq \eta(c,\delta)$$
\endproclaim
\demo{Proof}
The proof is based on the same idea as the proof of Lemma 7 in \cite{3}.
By a possible change of coordinates in the $x$ axis, we may and will
assume that $p(x)=\int_0^1 p(x,y)dy=1$, and 
that $|p(y|x)-1|\leq \delta'=2\delta/(1-\delta)$. 
Let $P_i$ be the law on $[0,1]$
with density $p(y|X_{i})$. 
Note that $P_i$ may be written as a mixture of a uniform law
(with weight $(1 -\delta')$) and another law on $[0,1]$,
depending on $X_i$ and
denoted $q_i$,
that is $P_i(dy)=(1-\delta')\lambda_1(dy) +\delta'q_i(dy).$
 Thus, the sample
$((X_{1},Y_{1}),\ldots,(X_{n},Y_{n}))$
possesses the same law as
$\tilde Z_n=((X_{1},(m_1U_1+(1-m_1)W_1)),
\ldots,(X_{n},(m_n U_n+(1-m_n)W_n)))$,
where $\{U_i\}_{i=1}^n$ is a sequence of i.i.d. uniform random
variables, independent of the sequence $\{X_i\}_{i=1}^n$,
$\{m_i\}_{i=1}^n$ is a sequence of i.i.d. Bernoulli$(1-\delta')$
 random variables,
independent of the sequences $\{U_i\}_{i=1}^n$ and $\{X_i\}_{i=1}^n$,
and $\{W_i\}_{i=1}^n$ is a sequence of random
variables whose law depends on the sequence $\{X_i\}_{i=1}^n$.
Let $I$ denote the set of indices with $m_i=1$, and let
$N_n=\sum_{i=1}^n 1_{m_i=1}=|I|$ denote the number of
indices where a uniform random variable is chosen
in the mixture.
Note that one may find a $\delta''=\delta''(\delta)\ge\delta', \delta''(\delta)\to_{\delta\to 0} 0$
such that
$$\limsup_{n\to\infty} \frac1n\log P(N_n/n<1-\delta'')=-\Big\{
(1-\delta'')\log\frac{1-\delta''}{1-\delta'}+\delta"
\log\frac{\delta''}{\delta'}\Big\}
<-3H_0(c)\,\tag{2.5}$$
for all $\delta$ small enough. Let $\tilde \ell_{\max}(n)$ denote the length
of the maximal increasing subsequence corresponding
to $\tilde Z_n$, then $\tilde \ell_{\max}(n)$ possesses the same
law as $\ell_{\max}(n)$ and, on the other hand, is not smaller
than the length of the maximal increasing subsequence when one considers
only those indices $i\in I$.  The latter is distributed precisely as
the length of the maximal increasing subsequence of a {\it uniform}
sample of random length $N_n$ which is independent of the
uniform sequence.
Therefore,
$$ P(\ell_{\max}(n)<(2+c)\sqrt{n})\leq
P(L_{\max}(n(1-\delta''))<(2+c)\sqrt{n})
+P(N_n/n<1-\delta'')\,.\tag{2.6}$$
The continuity of $H_0(c)$ implies that for $\delta$ small enough,
$$2H_0(\frac{2+c}{\sqrt{1-\delta''}}-2)<3H_0(c)\,.$$
Hence, (2.5), (2.6) and Theorem 1 imply that for $\delta$ small enough,
$$\limsup_{n\to\infty}  P(\ell_{\max}(n)<(2+c)\sqrt{n})\leq
-2H_0(\frac{2+c}{\sqrt{1-\delta''}}-2)=
-2H_0(c)+g(c,\delta)\,,$$
where the continuity of $H_0(\cdot)$ implies the required properties 
of $g(c,\delta)$.
The complementary lower bound is proved by a similar coupling.
\qed \enddemo

We now turn to general case and prove first a lower bound estimate:
for fixed $d>0 $ set
$$I _\mu(d )\equiv\inf\{ H(\nu|\mu): \nu\in\Cal M_1(\bold Q),
 2\bar J_\nu= d\}\,,\tag{2.7}$$
where $\Cal M_1(\bold Q)$ is the set of probability measures on $\bold Q$ and
$ H(\nu|\mu)$ denotes the relative entropy of $\nu$ with respect to $\mu$:
$$ H(\nu|\mu)=\int_{\bold Q}\log\frac{q(x,y)}{p(x,y)} \nu(dx,dy)$$
if $\frac{d\nu}{d\lambda}=q$ and $ H(\nu|\mu)=\infty$ otherwise.

Although an explicit computation for $I_\mu$ seems impossible,
it is quite easy to verify that
$I_\mu(d )=0$ for $d\ge 2\bar J_\mu$,
 and $0<I_\mu(d )<\infty$ if $0<d<2\bar J_\mu$,
(e.g., by combining Lemma 1 and Proposition 2.2 below).
Note that (1.1) implies that  
under $Q\equiv
\prod \nu$,
 for each $\epsilon>0$
$$\lim_{n\to\infty}Q(\ell_{\max}(n)\le (2\bar J_\nu+\epsilon)\sqrt{n}
)=1.$$
Using a standard change of measure argument, we get from this:
\proclaim{Proposition 2.1} For fixed $0<d<2\bar J_\mu$,
$$\liminf_{n\to\infty}\frac1n\log P(\ell_{\max}(n)\le d\sqrt n)\ge - I_\mu(d).$$
\endproclaim
However, a  simple
comparison with $2H_0$ in case
 $\mu=\lambda$ shows that $I_\lambda(2-\cdot)$ is not the correct rate function:
\proclaim{Lemma 1} Take $\mu=\lambda$, then
$$
\liminf_{\delta\searrow 0}
\frac{I_\lambda(2-\delta)}{\delta^3}\ge \frac{4}{9}\tag{2.8}$$
\endproclaim
\demo{Proof}
Assume the existence of  $\nu_\delta$ such that $2\bar J_{\nu_\delta}<2-\delta$
but $\lim_{\delta\to 0} H(\nu_\delta|\lambda)/\delta^3 <4/9$.
 For a fixed $K>0$ (independent of $\delta$),
let  $q_\delta(x,y)=d\nu_\delta/d\lambda(x,y)$, and denote
$$A_K=\{ (x,y)\in {\bold Q}: q_\delta(x,y) \leq (1+\delta K)\}.$$
One easily checks that $\lambda(A_K^\complement)<\delta/K^2$. Thus,  we
may assume
that  $q_\delta(x,y)=1+\delta m(x,y)$ for some
$m$ which, on $A_K$, is  bounded above by $K$.
Consider the set of
curves $x:(0,1-y)\longrightarrow\phi_y(x)=y+x$
where $0<y<\delta/3.$
Then
$$\align
1-\frac\delta2&>J_{\nu}(\phi_y) = \int_0^{1-y}\sqrt{1+\delta m(x,y+x)}\,dx \\
&\geq (1-y)+\frac{\delta}{2}\int_0
^{1-y}(K\wedge m(x,y+x))\,dx+O(\delta^2)
\endalign$$
and therefore
$$\align\frac{\delta}{3}-\frac16\delta^2&>
\int_0^{\delta/3} J_{\nu}(\phi_y)\,dy\\
&\geq \frac{\delta}{3}-\frac{1}{18}\delta^2
+\frac{\delta}{2}\int_0^{\delta/3}
\int_0^{1-y}(K\wedge m(x,y+x))\,dx\,dy +O(\delta^3)\,.\endalign
$$
Thus
$$\delta\int_0^{\delta/3}
\int_0^{1-y}(K\wedge m(x,y+x))\,dx\,dy<-2\delta^2/9+O(\delta^3)\,,$$
and, by symmetry,
for  $\Delta_{\delta}\equiv\{(x,y)\in
\bold Q: -\delta/3<x-y<\delta/3\}\cap A_K$,
 $$\nu(\Delta_{\delta})<
\lambda(\Delta_{\delta})-4\delta^2/9+O(\delta^3)\,,$$
with $\lim_{K\to\infty} \lim_{\delta\to 0} \lambda(\Delta_\delta)/\delta=2/3$.
Now the infimum of $H(\nu|\lambda)$
under the above condition is achieved at the constant density 
$1-\delta'$, where
$$\delta'=\frac{4\delta^2/9+
O(\delta^3)}
{\lambda(\Delta_\delta)}=\delta(2/3+g_K))
+O(\delta^2)$$
on $\Delta_{\delta}$, and $g_K\to_{K\to\infty} 0$
is a constant independent of $\delta$ whose value may change from line to line.
Substituting in $H(\nu|\lambda)$, one obtains
$$H(\nu_{\delta}|\lambda)\ge 
4\delta^3/9 +g_K \delta^3 + O(\delta^4)\,.$$
Taking the limits as $\delta\to 0$ (first) and then $K\to \infty$ yields a contradiction.
\qed\enddemo

Note that our argument is quite
rough and with additional work one could possibly identify
the constant $b\in[\frac49,\frac32]$
such that $\lim_{\delta\searrow 0}\frac{I_\lambda(2-\delta)}{\delta^3}=b$,
but this is quite irrelevant since a simple computation shows
$2H_0(0)=2H'_0(0)=2H^{''}_0(0)=0$ and $2H^{'''}_0(0)=\frac{1}{2}$ and
therefore
$$\lim_{\delta\uparrow 0}
\frac{2H_0(\delta)}{\delta^3}=\frac{1}{12}<\frac49\le\liminf_{\delta\searrow 0}
\frac{I_{\lambda}(2-\delta)}{\delta^3}.$$
Our next result shows a volume upper bound:
\proclaim{Proposition 2.2} Let $c<0$, then
$$\limsup_{n\to\infty}\frac1n\log P(\ell_{\max}(n)\le(2\bar J_\mu
+c)\sqrt n)<0.$$
\endproclaim
\demo{Proof}
Let $\phi$ denote an optimizer in (1.2) (whose existence is
ensured by \cite{3}). Fix $\Delta>0$ with $\Delta^{-1}$
an integer, and for $i=1,\ldots,\Delta^{-1}$
let 
$$\bold Q_i=[(i-1)\Delta,i\Delta]\times[\phi((i-1)\Delta),\phi(i\Delta)],\quad\rho_i=\Delta (\phi(i\Delta)-\phi((i-1)\Delta))p(i\Delta,\phi(i\Delta)).$$
Set
$n_i=n\rho_i$, 
and let $m_i$ denote the (random) number
 of points in the sample $\{Z_i\}_{i=1}^n$ which belong to $\bold Q_i$.
Then, for any $\epsilon>0$, by Sanov's theorem,
{ for $\Delta>0$ small enough},
$$\limsup_{n\to\infty}\frac1n\log P(\cup_{i=1}^{\Delta^{-1}}
A_i(\epsilon))
<0\quad\text{where}\quad A_i(\epsilon)=\{|n_i-m_i|>\epsilon n_i\},.\tag{2.10}$$
Let $\ell_i(m_i)$ denote the length of the longest increasing subsequence
among the $m_i$ points in $\bold Q_i$.
Let $P_{m_i}$ denote the law of the sample in $\bold Q_i$,
conditioned on $m_i$. Then $P_{m_i}$ possesses a product law with density
$p_i$, satisfying $ \lim_{\Delta\to 0} \sup_{x,y\in Q_i}|p_{i}(x,y)-1|=0$, c.f.
Lemma 2  of \cite{3}.
Corollary 1 implies that
for any $\epsilon'>0$ and all $i$, and all $\Delta$ small enough,
$$\limsup_{m_i\to\infty} \frac{1}{m_i} \log P_{m_i}(\ell_i(m_i)<2(1-\epsilon')
\sqrt{m_i})<0.$$
Hence, for $\Delta$ small, on $\cap_i\{|n_i-m_i|\le\epsilon n_i\}$,
$$\limsup_{n\to\infty}\frac1n\log P(\ell_i(m_i)<2(1-\epsilon')
\sqrt{m_i}\,|\{m_i\}_{i=1}^{\Delta^{-1}})<0\,.
\tag{2.11}$$
Note that  
(c.f. \cite{3}) $\sum_i\sqrt{\rho_i}\to_{\Delta\to 0} \bar J_\mu
$.
Choose $\Delta$  small enough such that
$|\sum_i\sqrt{\rho_i}-\bar J_\mu|<|c| /2$. Then,
$$\align
 P(\ell_{\max}(n)&\le (2\bar J_\mu+c)\sqrt{n}) \le
P(\cup_{i=1}^{\Delta^{-1}}A_i(\epsilon))\\
& \;\;\;+E[\cup_{i=1}^{\Delta^{-1}}P(\ell_i(m_i)<2\sqrt{n_i}+c\Delta \sqrt{n}/2\,|\{m_j\});\cap_{i=1}^{\Delta^{-1}}A_i(\epsilon)^\complement]\\
&\leq \Delta^{-1} \max_i P(|n_i-m_i|>\epsilon n_i)\\
&\;\;\;+ \Delta^{-1} \max_i\max_{m_i: |m_i-n_i|\le\epsilon n_i} P_{m_i}(\ell_i(m_i)<2\sqrt{m_i}(1+\epsilon)
+c'\sqrt{m_i})\,,\endalign $$
where $c'<0$ is independent of $\epsilon$ and $\Delta$.
Choosing now $\epsilon$ small enough such that $2\epsilon+c'<0$
and using (2.10) and (2.11), the proposition follows.
\qed
\enddemo

\noindent
{\bf Remark:} It is instructive  to relate $\ell_{\max}(n)$ to $\bar J_\nu$ for an measure $\nu$ associated with $R_n\equiv\frac1n\sum_{i=1}^n\delta_{Z_i},$ the empirical measure of the sample. 
 To this end, define for $\epsilon>0$, the random measure 
 $R_{n,\epsilon}$  with  constant density $\epsilon^{-2}$ on 
the squares
$\bold Q_\epsilon(Z_i)=
[-\frac\epsilon2,\frac\epsilon2]^2+Z_i,\ i=1,...,n,$
that is
$$ \frac{d\,R_{n,\epsilon}}{d\,\lambda}(x,y)=
\frac1{n}\sum_{i=1}^n
\epsilon^{-2}{\bold 1}_{\bold Q_\epsilon(Z_i)}(x,y)\,.
$$
Note that, $P$
almost surely $\epsilon_n\equiv\frac12\min_{1\le i<j\le n}(|X_i-X_j|\wedge |Y_i-Y_j|)>0.$ A simple computation shows
$$\frac{\ell_{\max}(n)}{\sqrt n}=\bar J_{R_{n,\epsilon_n}},\tag{2.12}$$
and therefore
$\{\ell_{\max}(n)\le \sqrt{n} d\}=\{\bar J_{R_{n,\epsilon_n}}\le d\},$
for each $d>0$.
However a derivation of  the large deviation principle using
this equality fails, due to the discontinuity of the mapping
$\nu\longrightarrow\bar J_\nu$. In particular,
 $R_{n,\epsilon_n}$ converges weakly to $\mu$, on the other hand, we have $\lim_{n\to\infty}\frac{\ell_{\max}(n)}{\sqrt n}=2\bar J_\mu.$

\heading{\S. 3 The upper tail}\endheading

Here the situation is quite different from the lower tail,
and in some sense much simpler. Our
first result is:
\proclaim{Theorem 2} For all $c\ge0$
$$\lim_{n\to\infty}\frac1{n^{1/2}}\log P(L_{\max}(n)\ge( 2+c)\sqrt n)=
-U_0( c),\tag{3.1}$$
where $U_0:[0,\infty)\longrightarrow [0,\infty)$ is a continuous,
strictly increasing convex function
with $U_0( c)=0$ iff $c=0$,
and 
$$ U_0(c)=\beta(c):= 2(2+c)\cosh^{-1}(c/2+1)-2
\sqrt{c^2+4c}\,.\tag{3.2}$$
\endproclaim

\noindent
{ Note that $U_0(c)=O(c^{3/2})$ as $c\to 0$, this is also  predicted by
the behavior $2H_0(c)= O(c^3)$ as $c\to 0$ . The explicit computation
of $U_0(c)$ was first done in \cite{8} using Hammersley's
particle system.}
\demo{Proof}
As pointed out by \cite{1},
the convergence  in (3.1),
and the convexity and monotonicity
of $U_0$ follows from sub-additivity.
We briefly recall the argument.
Let ${\Cal N}_n$ denote the number of points in a Poisson 
point process of rate $\lambda_n=n\lambda$ on the unit square, and let
$\bar L_{\max}(\Cal N_n)$ denote the longest increasing subsequence
in that sample. Then, for any $\epsilon>0$, a direct computation
using the Poisson distribution yields
$$\limsup_{n\to \infty} \frac1{\sqrt{n}}
\log P(|{\Cal N}_n/n-1|>\epsilon)=-\infty\,.\tag{3.3}$$

On the other hand, conditioned on ${\Cal N}_n$, the law of the sample is
uniform and hence
$$P(\bar L_{\max}({\Cal N}_n)=x|{\Cal N}_n=m)=P(L_{\max}(m)=x)\,.$$
Therefore,
$$ P(\bar L_{\max}({\Cal N}_n)\geq(2+c)\sqrt{n})/
P({\Cal N}_n>n(1-\epsilon))\geq P(L_{\max}(n(1-\epsilon))>(2+c)\sqrt{n})$$
while
$$ P(\bar L_{\max}({\Cal N}_n)\geq(2+c)\sqrt{n})
-P({\Cal N}_n>n(1+\epsilon))\leq P(L_{\max}(n(1+\epsilon))>(2+c)\sqrt{n})\,,$$
which implies (using (3.3))
that (3.1) holds as soon as it
holds with $\bar L_{\max}(\Cal N_n)$ replacing $L_{\max}(n)$.

On the other hand, consider the squares
$\bold Q^1=[0,\sqrt{n}/(\sqrt{n}+\sqrt{m})]^2$
and 
$\bold Q^2=(\sqrt{n}/(\sqrt{n}+\sqrt{m}),1]^2\subseteq\bold Q$, 
and denote by $\bar L^1_{\max}(\Cal N_{(\sqrt n +\sqrt m)^2}) $ and
$\bar L^2_{\max}(\Cal N_{(\sqrt n +\sqrt m)^2})$ 
 the length
of the longest increasing subsequence in the squares
$\bold Q^1$ and $\bold Q^2$, 
corresponding to $\Cal N_{(\sqrt n+\sqrt m)^2}$.
The scaling and independence properties of the Poisson process imply that
$\bar L^i_{\max}(\Cal N_{(\sqrt n +\sqrt m)^2})$, $i=1,2$, 
are independent, and that the laws of 
$\bar L^1_{\max}(\Cal N_{(\sqrt n +\sqrt m)^2}) $
and $L_{\max}({\Cal N}_n)$,  respectively 
$\bar L^2_{\max}(\Cal N_{(\sqrt n +\sqrt m)^2}) $ and $L_{\max}({\Cal N}_m)$, 
are identical. Therefore,
since 
$$\bar L_{\max}({\Cal N}_{(\sqrt{n}+\sqrt{m})^2})
\geq\bar L^1_{\max}(\Cal N_{(\sqrt n +\sqrt m)^2}) +\bar L^2_{\max}(\Cal N_{(\sqrt n +\sqrt m)^2})\,,$$
we deduce that
$$\align
&P(\bar L_{\max}({\Cal N}_{(\sqrt{n}+\sqrt{m})^2})>(2+c)(\sqrt{n}+\sqrt{m}))
 \\
&\;\;\;\geq 
P(\bar L^1_{\max}(\Cal N_{(\sqrt n +\sqrt m)^2}) +
\bar L^2_{\max}(\Cal N_{(\sqrt n +\sqrt m)^2}) 
>(2+c)(\sqrt{n}+\sqrt{m}))
\\
&\;\;\;
\geq P(\bar L^1_{\max}(\Cal N_{(\sqrt n +\sqrt m)^2}) >(2+c)\sqrt{n})
P(\bar L^2_{\max}(\Cal N_{(\sqrt n +\sqrt m)^2}) >(2+c)\sqrt{m})\,,\\
&\;\;\;
= P(\bar L_{\max}({\Cal N}_n)>(2+c)\sqrt{n})
P(\bar L_{\max}({\Cal N}_m)>(2+c)\sqrt{m})\,,\endalign $$
which immediately implies the existence and convexity of the limit
$$\bar U_0(c)=\lim_{k\to\infty} \frac1k\log P(\bar L_{\max}({\Cal N}_{k^2}
)>(2+c)k)\,.$$
Next, since 
$$\align 
P(\bar L_{\max}(\Cal N_{([\sqrt{n}]+1)^2}>
(2+c)\sqrt{n})& \geq
P(\bar L_{\max}(\Cal N_n)>(2+c)\sqrt{n})\\
& \geq P(\bar L_{\max}(\Cal N_{[\sqrt{n}]^2}>
(2+c)\sqrt{n})\,,\endalign $$
(3.1) follows with $U_0=\bar U_0$. It thus remains only to 
explicitely compute $U_0(c)$.

In fact, Kim has already observed that $U_0(c)\geq \beta(c)$,
$c\geq 0$, see equation (1.6) in \cite{5}. We thus concentrate in the
sequel in the reverse inequality. The proof is constructive:
we exhibit an appropriate collection of Young shapes.

Fix $n$ large enough, $\epsilon>0$ small
and $c>0$, and let ${\Cal T}_{c,n}$ denote 
the set consisting of Young shapes of size $n_c=n-\lceil(c+\epsilon)
\sqrt{n}\rceil$. Recall the function
$f_0^0\in \Cal F $ defined in Section 2, and define
$$ {\Cal T}_{c,n}^\epsilon=\{\tau\in {\Cal T}_{c,n}:
|\tau(0)-2\sqrt{n}|\leq \epsilon\sqrt{n}, \sup_x|
n_c^{-1/2}\tau(x\sqrt{n_c})-f_0^0(x)|
\leq
\epsilon\}\,.$$
It follows easily from \cite{5}, \cite{6} and (2.3) that
for $n$ large enough,
$$ \sum_{\tau\in {\Cal T}_{c,n}^\epsilon} 
\frac{n_c!}{(\pi(\tau))^2}\geq \frac12\,. \tag{3.4}$$
For each $\tau\in {\Cal T}_{c,n}^\epsilon$, define a new  Young shape
$\tau'$ obtained by increasing the height  $\tau(0)$ by
$\lceil(c+\epsilon)\sqrt{n}\rceil$. Note that $\tau'(0)\geq (2+c)\sqrt{n}$
while $|\tau'|=n$. From (2.3), we have
$$ \align
P(L_{\max}(n)\geq (2+c)\sqrt{n})
&\geq \sum_{\tau':\tau\in {\Cal T}_{c,n}^\epsilon} 
\frac{n!}{(\pi(\tau'))^2}
\\
&=
\frac{n!}{n_c!}\sum_{\tau\in {\Cal T}_{c,n}^\epsilon}
\frac{n_c!}{(\pi(\tau))^2}\Big(\frac{\pi(\tau)}{\pi(\tau')}\Big)^2\,.
 \endalign
$$
For $\tau\in {\Cal T}_{c,n}^\epsilon$, 
$$\frac{\pi(\tau')}{\pi(\tau)} =
\Big(\prod_{i=1}^{\lceil(c+\epsilon)\sqrt{n}\rceil} i\Big)
\prod_{j=1}^{\tau(0)}
\frac{(\tau(0)+\lceil(c+\epsilon)\sqrt{n}\rceil-j+\tau(j))}
{(\tau(0)-j+\tau(j))}\,.$$
Note that, due to the definition of ${\Cal T}_{c,n}^\epsilon$,
$$
\align &\log \prod_{j=1}^{\tau(0)}
\frac{(\tau(0)+\lceil(c+\epsilon)\sqrt{n}\rceil-j+\tau(j))}
{(\tau(0)-j+\tau(j))}
\\
&\leq  \sqrt{n} 
\Big[
\int_0^2 \log(2+c-x+f^0_0(x))dx-\int_0^2 \log(2-x+f_0^0(x)) dx
+C_\epsilon\Big] +o(\sqrt{n})\,, \endalign$$
where $C_\epsilon\to_{\epsilon\to 0} 0$ does not depend on $n$.
Using the change of variables $x-f_0^0(x)=\xi, x=h_0(\xi)$ and
$g_0(\xi)=h_0(\xi)-\xi(1+\frak{sign}(\xi))/2$ as in Pg. 212 of \cite{6},
one obtains after some manipulations that 
$$ \align
&\int_0^2 \log(2+c-x+f^0_0(x))dx-\int_0^2 \log(2-x+f_0^0(x)) dx\\
& =\int_{-\infty}^\infty  g_0'(\xi)
(\log(2+c-\xi)-\log(2-\xi)) d\xi+
\int_0^2 (\log(2+c-\xi)-\log(2-\xi)) d\xi\\
&=\pi \tilde g_0(2+c)-\pi \tilde g_0(2)+
(2+c)\log(2+c)
-c\log c-2\log 2\,, \tag{3.5}\endalign$$
where $\tilde g_0$ denotes the Hilbert transform of $g_0$ and
is given by (2.31) in \cite{6}. Note however that by (2.22) 
in \cite{6}, $\pi\tilde g_0(2)=2-2\log 2$,
while (2.22) and (2.32) in \cite{6} imply
$\pi\tilde g_0(c+2)=(2+c)-(2+c)\log(2+c) + \beta(c)/2\,.$
Substituting in (3.5), and then using (3.4),
one concludes that
$$ \liminf_{n\to\infty} \frac{1}{\sqrt{n}} \log P(L_{\max}(n)\geq
(2+c)\sqrt{n})
\geq -(\beta(c)+2C_\epsilon)\,.$$
Taking $\epsilon\to 0$ yields the desired conclusion that $U_0(c)\leq
\beta(c)$
for $c>0$.
\qed
\enddemo

The following corollary follows from Theorem 2 in the same way
that Corollary 1 followed from Theorem 1:
\proclaim{Corollary 2} For any $c>0$ there exists a function
$\bar \eta(c,\delta)$ satisfying 
$$\lim_{\delta\to 0} \bar \eta(c,\delta)=0$$
such that if $p(x,y)$ satisfies $(1-\delta)\leq p(x,y)\leq (1+\delta)$ 
then
$$\limsup_{n\to\infty}
|\frac1{\sqrt{n}}\log P( \ell_{\max}(n)>(2+c)\sqrt{n})+U_0(c )|
\leq \bar \eta(c,\delta)\,.$$
\endproclaim

Let $K\subseteq B^\uparrow$ be the
set of solution to the variational problem (1.2).
\proclaim{Theorem 3} For all $c\ge 0$
$$\lim_{n\to\infty}\frac1{n^{1/2}}\log P(\ell_{\max}(n)
\ge (2\bar J_\mu+c)\sqrt n)=-\bar J_\mu U_0\big(
c/\bar J_\mu)\,.\tag{3.6}$$
Next assume that $K=\{\phi_1,...,\phi_r\}$.
Then for each $\delta>0$ and longest increasing
subsequence $Z^{\max}=\{(X_{i_j},Y_{i_j}), j=1,...,\ell_{\max}(n)\}\,,$
$$\lim_{n\to\infty}
P(\min_{\alpha=1}^r \max_{j=1}^{\ell_{\max}(n)} |Y_{i_j}-\phi_\alpha(X_{
i_j})|<\delta\big |  \ell_{\max}(n)\ge (2\bar J_\mu+c)\sqrt n)=1.$$
\endproclaim
\demo{Proof}
We begin by 
providing a lower bound in (3.6). Let $\phi$ denote a maximizer
in (1.2), and define $\Delta,\rho_i,n_i,m_i,\bold Q_i,\ell_i(m_i)$ be as in
the beginning of the proof of Proposition 2.2. Fix $\delta>0$, and
reduce $\Delta$ if necessary. By
Sanov's theorem,
$$
\limsup_{n\to\infty} 
\frac{1}{n^{1/2}} \log P(|n_i-m_i|>\delta n_i)=-\infty\,.
$$
Hence,
$$ \align 
&\liminf_{n\to\infty} 
\frac1{n^{1/2}}\log P(\ell_{\max}(n)
\ge (2\bar J_\mu+c)\sqrt{ n})
 \geq  \\
&\liminf_{n\to\infty} 
\frac1{n^{1/2}}\log P(\sum_{i=1}^{\Delta^{-1}} \ell_i(m_i)
\ge (2\bar J_\mu+c)\sqrt{ n}; \cap_{i=1}^{\Delta^{-1}} \{
|n_i-m_i|\leq \delta n_i\})\,.\endalign$$
Next, for each $i$ and $t_i>0$, 
using $\lim_{\Delta\to 0} \sup_{x,y\in Q_i}|p_{i}(x,y)-1|=0$,
one has by Corollary 2 that
for some $\delta'(\Delta)\to_{\Delta\to 0} 0$, 
$$|\lim_{m_i\to\infty} \frac{1}{m_i^{1/2}} \log P_{m_i}
(\ell_i(m_i)>(2+t_i)\sqrt{m_i})
+U_0(t_i)|\leq\bar{\eta}(t_i,\delta') \,.\tag{3.7} 
$$
($P_{m_i}$ is
the law of the sample in $\bold Q_i$ conditioned on $m_i$).

Recall
$\rho_i=n_i/n=\Delta (\phi(i\Delta)-\phi((i-1)\Delta))
p(i\Delta,\phi(i\Delta))$,
and fix a sequence $\{t_i\geq 0\}_{i=1}^{\Delta^{-1}}$
such that
$$\sum_{i=1}^{\Delta^{-1}} 
(2+t_i)\sqrt{(1-\delta)\rho_i}\geq (2\bar J_\mu+c)\,.\tag{3.8}$$
Then, 
$$ \align 
&
\frac1{n^{1/2}}\log P(\sum_{i=1}^{\Delta^{-1}} \ell_i(m_i)
\ge (2\bar J_\mu+c)\sqrt{ n}; \cap_{i=1}^{\Delta^{-1}} \{
|n_i-m_i|\leq \delta n_i\})\\
\geq &
\frac1{n^{1/2}}\log \inf_{\{m_i: |n_i-m_i|\leq \delta n_i\}}
P(\sum_{i=1}^{\Delta^{-1}} \ell_i(m_i)
\ge (2\bar J_\mu+c)\sqrt{ n}| \{m_j\})\\
\ge &
\frac1{n^{1/2}}\log 
\prod_{i=1}^{\Delta^{-1}}
P_{n_i(1-\delta)}( \ell_i(m_i)
\ge (2+t_i)
\sqrt{(1-\delta)n_i})\endalign
$$
where the last inequality is a consequence of (3.8),
the monotonicity
of $\ell_i(m_i)$ in $m_i$, and of the (conditional in $\{m_j\}$)
independence of the $\ell_i(m_i)$. Thus,
combining the above with (3.7),
one concludes that
$$ 
\liminf_{n\to\infty} 
\frac1{n^{1/2}}\log P(\ell_{\max}(n)
\ge (2\bar J_\mu+c)\sqrt{ n})
\geq -\sum_{i=1}^{\Delta^{-1}} \sqrt{\rho_i (1-\delta)}
(U_0(t_i)+\bar \eta(t_i,\delta')\}  \,.$$
Since the last bound is valid for any choice of $\{t_i\}$ satisfying
(3.8), we conclude that
$$ \align 
&\liminf_{n\to\infty} 
\frac1{n^{1/2}}\log P(\ell_{\max}(n)
\ge (2\bar J_\mu+c)\sqrt{ n})
 \geq  \\
&-\inf\{\sum_{i=1}^{\Delta^{-1}} \sqrt{\rho_i}
(U_0(t_i)+\bar\eta (t_i,\delta')):
t_{\cdot}\geq 0, \sum_{i=1}^{\Delta^{-1}} (2+t_i)\sqrt{(1-\delta)
\rho_i}
\geq (2\bar J_\mu+c)\}\,. \tag{3.9}\endalign $$

Recall (c.f. \cite{3}) that
$\sum_{i=1}^{\Delta^{-1}}\sqrt{\rho_i}\to_{\Delta\to 0} \bar J_\mu$.
The smoothness
of $\phi$ proved in \cite{3} and (3.9) imply therefore, by
taking the limit $\Delta\to 0$ in the right hand side of (3.9),  that
$$\align
& \liminf_{n\to\infty} 
\frac1{n^{1/2}}\log P(\ell_{\max}(n)
\ge (2\bar J_\mu+c)\sqrt n)
\\
&\geq 
-\inf\{\int_0^1 \sqrt{\dot{\phi}(x)p(x,\phi(x))}U_0(t(x)) dx:\\
&\;\;\;t\in C([0,1],{\Bbb R^+})\,,
\int_0^1 t(x) \sqrt{\dot{\phi}(x) p(x,\phi(x))} dx\geq c\}\,.\tag{3.10}
\endalign $$
Making the change of variables $dy=\sqrt{\dot{\phi}(x) p(x,\phi(x))}dx,
t(x)\to \bar t(y)$, the right hand side of (3.10) becomes
$$-\inf\{ \int_0^{\bar J_\mu} U_0(\bar t(y)) dy\,:\,
\bar t\in C([0,\bar J_\mu],\Bbb R^+),
\int_0^{\bar J_\mu} \bar t(y) dy \geq c\}\,.\tag{3.11}$$
Take now $\bar t(y)=c/\bar J_\mu$ to conclude from (3.10)
that
$$ \liminf_{n\to\infty} 
\frac1{n^{1/2}}\log P(\ell_{\max}(n)
\ge (2\bar J_\mu+c)\sqrt n)
\geq - \bar J_\mu U_0(\frac{c}{J_\mu})\,.$$

The proof of the complimentary upper bound is only 
slightly more complicated, and
involves the same tools as in \cite{3}.
Let $\Delta_y<<\Delta$, with $\Delta_y^{-1}$ an integer.
Define
a ``block curve" as an integer valued
sequence $\{j(i)\}_{i=1}^{\Delta^{-1}}$,
satisfying $j(i+1)>j(i)$ and $j(\Delta^{-1}) \Delta_y\leq 1$.
Let ${\Cal B}^\Delta$ denote the set
of all possible block curves, and note
that the cardinality of ${\Cal B}^\Delta$ is finite.
To any block curve $b\in {\Cal B}^{\Delta}$ 
associate naturally a (piecewise linear) curve $\phi_b$,
and define
$\bar{\bold Q}_i=[(i-1)\Delta, i\Delta]\times[j(i-1)\Delta_y,(j(i)+1)\Delta_y]$,
$\bar m_i$ as the number of points within $\bar{\bold Q}_i$,
$\bar \rho_i=\Delta (j(i)-j(i-1)+1)\Delta_y p(i\Delta, j(i)\Delta_y)$, $\bar n_i=n
\bar\rho_i$ and
$\bar \ell_i(b,n)$ as the length of the longest increasing
subsequence within $\bar{\bold Q}_i$. Clearly,
$\ell_{\max}(n)\leq \max_{b\in {\Cal B}^\Delta} \sum_{i=1}^{\Delta^{-1}}
\bar \ell_i(b,n)$. Hence,
$$\align
&\frac1{n^{1/2}}\log P(\ell_{\max}(n)
\ge (2\bar J_\mu+c)\sqrt n)
\leq\\
& \frac1{n^{1/2}}\log |{\Cal B}^\Delta|+
\max_{b\in {\Cal B}^\Delta}
\frac1{n^{1/2}}\log P(\sum_{i=1}^{\Delta^{-1}} \bar \ell_i(b,n)
\geq (2\bar J_\mu+c)\sqrt n)\,.\endalign$$
Fix $\delta>0$ small. Repeating the argument used
in the proof of the lower bound, one finds (reducing $\Delta, 
\Delta_y/\Delta$ if necessary, but independently of $n$) that
$$\align 
&\limsup_{n\to\infty}
\frac1{n^{1/2}}\log P(\sum_{i=1}^{\Delta^{-1}} \bar \ell_i(b,n)
\geq (2\bar J_\mu+c)\sqrt{n})\\
&
\leq - \inf\{\sum_{i=1}^{\Delta^{-1}} \sqrt{\bar \rho_i} (U_0(t_i)+
\bar \eta(\delta',t_i))\,:\,
t_i\geq 0, \sum_{i=1}^{\Delta^{-1}} (2+t_i))
\sqrt{\bar \rho_i(1+\delta)}\geq 2\bar J_\mu+c\}\,.\endalign $$

Let $J_\phi=\int_0^1 \sqrt{\dot{\phi}(x)p(x,\phi(x))} dx$.
With $\Delta$ small enough,
$$|\sum_{i=1}^{\Delta^{-1}}\sqrt{\bar \rho_i}- J_{\phi_b}|
\leq \delta\,.$$
Hence, taking now first $n\to\infty$ and then $\Delta\to 0$, 
followed by $\delta\to 0$, one concludes that
$$\align
& \limsup_{n\to\infty} 
\frac1{n^{1/2}}\log P(
\ell_{\max}(n)
\geq (2\bar J_\mu+c)\sqrt n)\\
& \leq
-\inf_{\phi\in B^\uparrow} \inf \{\int_0^1 \sqrt{\dot{\phi}(x) p(x,\phi(x))}
U_0(t(x))dx\,:\,\\
& \;\;\;t(\cdot)\geq 0, 
\int_0^1 t(x)\sqrt{\dot{\phi}(x) p(x,\phi(x))} dx
\geq 2(\bar J_\mu-J_\phi)+c)\,.\tag{3.12}\endalign $$
Making the same change of variables as in the proof of the lower bound, 
the right hand side of (3.12) equals
$$
\align
I&=\inf_{\phi\in B^\uparrow} \inf\{\int_0^{J_\phi} U_0(\bar t(y))dy\,:\,
\int_0^{J_\phi} \bar t(y) dy\geq c+2(\bar J_\mu-J_\phi)\}\\
&\geq
\inf_{\phi\in B^\uparrow} 
\{J_\phi U_0( -2+\frac{c+2\bar J_\mu}{J_\phi})\}
\,,\tag{3.13}\endalign $$
where the last inequality follows from 
the convexity of $U_0$ and
Jensen's inequality.
Let $x= -2+\frac{c+2\bar J_\mu}{J_\phi}$. Then, since $J_\phi\leq \bar J_\mu$,
$x\geq c/\bar J_\mu$.
Hence, using again the convexity of $U_0$ and the fact that $U_0(0)=0$,
$cU_0(x)\geq \bar J_\mu xU_0(c/\bar J_\mu)$
and hence
$$J_\phi U_0(x)\geq\frac{ x J_\phi \bar J_\mu}{c} U_0(c/\bar J_\mu)
\geq \bar J_\mu U_0(c/\bar J_\mu)\,,\tag{3.14}$$
with the second inequality being strict unless $J_\phi=\bar J_\mu$.
(3.12), (3.13) and (3.14) imply the required upper bound.

Finally, the last statement of Theorem 3 follows from the fact that 
the inequality in (3.14) is strict unless $J_\phi=\bar J_{\mu}$,
and the fact that the assumption of finite $K$ implies A4 of \cite{3}
(the proof is similar to the proof of Lemma 4 in \cite{3} and is thus omitted).
\qed 
\enddemo

\heading{\S 4 Open problems and remarks}\endheading
We conclude this paper with a list of comments and open problems.

\bigpagebreak
\noindent
1) We have left open the question of
existence of limit in (1.4) and of the computation of $H_\mu(c)$.
After a discretization as used in Section 2,
maybe techniques borrowed from percolation may allow one to control
the interaction between overlapping   ``block curves". 

 \bigpagebreak

\noindent
2) We have seen in Theorem 3 that,
 under the conditioning $\{\ell_{\max}(n)\ge(2\bar J_\mu+c)\sqrt n\}, (c>0),$ 
any longest increasing subsequence concentrates along the solution to the variational problem (1.2). The corresponding question for
$\{\ell_{\max}(n)\le(2\bar J_\mu+c)\sqrt n\}, (-2\bar J_\mu<c<0),$ remains unsolved,
even in case $\mu=\lambda$. 

\bigpagebreak

\noindent
3) Under the assumption that $K$ is finite,
the strict convexity of $U_0(c)$ implies uniqueness of
 the minimizing function $t(\cdot)$ in
(3.10) and gives the profile of the longest increasing subsequence
under the conditioning that an upper  tail deviation occurred, for
any $\mu$. Indeed, for an optimal curve $\phi$,
the minimizing function $t(x)=t_\phi(x)$ in the variational
problem (3.10) is readily   seen to have the interpretation as
the (local) fluctuation from the mean behavior, and strict convexity
of $U_0$ would imply that $t(x)=c/\bar J_\mu$, a constant.
\comment
\noindent
{\it Note: after this work was completed, T. Sepp\"{a}l\"{a}inen
has kindly provided us with a copy of \cite{8}, which contains
an explicit computation of $U_0(\cdot)$. In particular, it follows
from \cite{8} that $U_0(\cdot)$ is strictly convex.}
\endcomment
\bigpagebreak

\noindent
4) It is natural to ask what happens when $\bold Q$ is replaced by
$[0,1]^d$, $d>2$.
The subadditivity argument for the upper tail is the
same, as well as the analog of Theorem 3 (with exponential speed
$n^{1/d}$, and functional $J_\mu$ as given in 
\cite{3}, page 864). What about the lower tail?  The lack of a 
direct probabilistic proof of  Theorem 1, and the unavailability of
the Schensted correspondence  in higher dimension makes
finding the analog of Theorem 1 challenging.
One can still show, however, that
the order of decay is exponential in $n$.
\bigpagebreak 

\noindent
5) As pointed out to us by P. Baxendale,
it seems reasonable to expect that under the conditioning
$\ell_{\max}(n)$ $\geq $ $\alpha_n$, $\alpha_n/$ $\sqrt{n} $ $\to \infty$,
the maximizing subsequences concentrate around  the solutions of 
the optimization problem (1.2).  For $\alpha_n=n$, this was
proved in \cite{3}, and the technique of the proof seems
to carry to the general case.

\bigpagebreak
\noindent
6) Let $\Cal N_n$ be the number of points of a Poisson point process
on $\bold Q$ with intensity $n\mu$ and denote by $\bar\ell_{\max}(\Cal N_n)$
the length of the longest increasing subsequence of $\Cal N_n$.
Note that, conditioned on $\Cal N_n=m$, the law of $\bar\ell_{\max}(\Cal N_n)$
is the same as the law of $\ell_{\max}(m)$. 
Applying the same type of argument as in the first step of the proof of 
Theorem 2, one shows 
$$\lim_{n\to\infty}\frac1{n^{1/2}}\log P(\bar\ell_{\max}(\Cal N_n)
\ge (2\bar J_\mu+c)\sqrt n)=-\bar J_\mu U_0\big(
c/\bar J_\mu).$$
The corresponding result for the lower tail (in the case $\mu=\lambda$)
can also be read off Theorem 1, c.f. \cite{8}. Note that in this case
the rate function does differ from the uniform case due to fluctuations
in the number of points in the Poisson sample.\

\end